\documentclass[aps,prl,reprint,superscriptaddress,showpacs]{revtex4-1}

\usepackage[pdftex]{graphicx}
\setkeys{Gin}{clip=true,draft=false}
\DeclareGraphicsExtensions{.pdf}
\usepackage[cmex10]{amsmath}
\usepackage{amsthm}
\usepackage{amsfonts}
\usepackage[cspex,bbgreekl]{mathbbol}
\newcommand{\R}{{\ifmmode\mathbb{R}\else$\mathbb{R}$\fi}}
\newcommand{\K}{{\ifmmode\mathbb{K}\else$\mathbb{K}$\fi}}
\newcommand{\N}{{\ifmmode\mathbb{N}\else$\mathbb{N}$\fi}}
\newcommand{\C}{{\ifmmode\mathbb{C}\else$\mathbb{C}$\fi}}
\newcommand{\D}{\text{d}}

\newcommand{\I}{\text{i}}

\newcommand{\gei}{{\I}}

\newcommand{\fracd}[2]{\frac{\displaystyle #1}{\displaystyle #2}}

\newcommand{\totder}[2]{\frac{\displaystyle\D #1}{\displaystyle\D #2}}

\newcommand{\abs}[1]{{\left| #1 \right|}}

\newcommand{\up}[1]{{\text{#1}}}
\newcommand{\idest}{{i.e.}}

\newcommand{\bvec}{\mathbf{b}}

\newcommand{\fvec}{\mathbf{f}}

\newcommand{\zvec}{\mathbf{z}}

\newcommand{\Bmat}{\mathbf{B}}
\newcommand{\Cmat}{\mathbf{C}}
\newcommand{\Dmat}{\mathbf{D}}

\newcommand{\Amat}{\mathbf{A}}
\newcommand{\Xmat}{\mathbf{X}}

\newcommand{\Fmat}{\mathbf{F}}

\newcommand{\xvec}{\mathbf{x}}
\newcommand{\gvec}{\mathbf{g}}
\newcommand{\vvec}{\mathbf{v}}
\newcommand{\yvec}{\mathbf{y}}
\newcommand{\evec}{\mathbf{e}}

\newcommand{\kvec}{\mathbf{k}}
\newcommand{\avec}{\mathbf{a}}
\newcommand{\zerovec}{\mathbf{0}}

\newcommand{\cvec}{\mathbf{c}}

\newcommand{\nvec}{\mathbf{n}}

\newcommand{\rvec}{\mathbf{r}}
\newcommand{\Rvec}{\mathbf{R}}

\newcommand{\Mmat}{\mathbf{M}}

\newcommand{\Qmat}{\mathbf{Q}}
\newcommand{\Rmat}{\mathbf{R}}

\newcommand{\Kmat}{\mathbf{K}}

\newcommand{\esp}{\operatorname{\mathrm{e}}}
\newcommand{\espo}[1]{\esp^{\displaystyle #1}}

\newcommand{\lambdavec}{\boldsymbol{\lambda}}

\DeclareMathOperator{\Lt}{\mathcal{L}_\mathcal{T}}

\DeclareMathOperator{\Ltp}{\mathcal{L}'_\mathcal{T}}
\DeclareMathOperator{\Lts}{\mathcal{L}''_\mathcal{T}}

\DeclareMathOperator{\T}{\mathcal{T}}
\DeclareMathOperator{\Tb}{\mathcal{T}}

\DeclareMathOperator{\HT}{\mathcal{H}}
\DeclareMathOperator{\V}{\mathcal{V}}
\DeclareMathOperator{\Vtilde}{\tilde{\mathcal{V}}}

\newtheorem*{theorem}{Theorem}

\begin{document}

%\preprint{}

\title{Generalized Floquet theory: application to dynamical systems with memory and Bloch's theorem for nonlocal potentials}


\author{Fabio L. Traversa}
\email[]{fabiolorenzo.traversa@uab.es}
%\thanks{}
\affiliation{Departament d'Enginyeria Electr\`onica, Universitat Aut\`onoma de Barcelona, 08193-Bellaterra (Barcelona), Spain}
\affiliation{Department of Physics, University of California, San Diego, La Jolla, California 92093, USA}

\author{Massimiliano Di Ventra}
\email[]{diventra@physics.ucsd.edu}
%\thanks{}
\affiliation{Department of Physics, University of California, San Diego, La Jolla, California 92093, USA}

\author{Fabrizio Bonani}
\email[]{fabrizio.bonani@polito.it}
%\thanks{}
\affiliation{Dipartimento di Elettronica e Telecomunicazioni, Politecnico di Torino, 10129 Torino, Italy}


\date{\today}

\begin{abstract}
Floquet theory is a powerful tool in the analysis of many physical phenomena, and extended to spatial coordinates provides the basis for Bloch's theorem. However, in its original formulation it is limited to linear systems with periodic coefficients. Here, we extend the theory by proving a theorem for the general class of systems including linear operators commuting with the period-shift operator. The present theorem greatly expands the range of applicability of Floquet theory to a multitude of phenomena that were previously inaccessible with this type of analysis, such as dynamical systems with memory. As an important extension, we also prove Bloch's theorem for nonlocal potentials.
\end{abstract}

\pacs{05.45.-a, 71.20.-b, 02.30.Oz, 02.30.Sa}
% 02.30.Oz 	Bifurcation theory
% 02.30.Sa 	Functional analysis
% 05.45.-a 	Nonlinear dynamics and chaos
% 84.30.Bv 	Circuit theory
%\keywords{}

\maketitle

Floquet theory is a century-old theory that allows an easy analysis of the solutions of linear differential equations that have periodic coefficients \cite{Floquet,Farkas}. These types of equations are very common in several areas of science and technology, and as such the applications of Floquet theory range from quantum \cite{Shirley,buttiker} to classical \cite{Gammaitoni,Haken} physics, chemistry \cite{Boland}, electronics \cite{Demir1,noiFN,noi}, dynamical systems \cite{Guckenheimer} and many more. However, its original formulation is limited to linear systems with periodic coefficients. Consequently, it is a powerful tool to study nonlinear perturbations, noise, and stability of systems depending {\it instantaneously} on time, and admitting periodic steady states.

With reference to condensed matter systems, the classical Bloch theorem \cite{bloch} can be derived as a simple
 corollary of Floquet theorem applied to space-dependent periodic potentials. Bloch states play a fundamental role in defining the concepts of allowed and forbidden energy bands in crystalline solids under the form of the crystal dispersion relation~\cite{Ashcroft-Mermin}. However, also in this case the usual Bloch theorem is applied to potentials depending \textit{locally} on position.

These limitations leave out many important phenomena that cannot be simply described using these types of equations.
This is the case, for example, of all dynamical systems whose state depends on their past dynamics. These systems are the norm rather than the exception in natural phenomena \cite{Max_review} and their applications span a wide range of problems (see, e.g., \cite{Max_Proc_IEEE} for recent examples). Turning to condensed matter physics, non-local potentials again naturally arise whenever the many body quantum problem is simplified by tracing out several atomic or electronic degrees of freedom -- e.g., those associated with the ``core'' electrons. \cite{Ashcroft-Mermin,MaxBook}.
It would be thus of great value to extend Floquet theory to such a general class of physical systems.

In this letter, we state and prove a generalization of Floquet theorem that applies to systems with memory.
As a corollary of such a theorem, we prove Bloch theorem for 3D crystals where a non-local potential energy profile is present, thus showing that the fundamental tool we provide is useful to probe the properties of a wide range of physical systems.

\paragraph{Introductory remarks.}  Let us first consider a general independent variable $\sigma$, that can represent either the time $t$ or a scalar spatial variable $x$, since the result we derive can be applied both to time- and (one-dimensional) space-dependent equations. The standard Floquet theorem applies to linear homogeneous systems of the form $\D \zvec_0/\D \sigma=\Amat(\sigma)\zvec_0(\sigma)$ where $\zvec_0: \R\rightarrow\R^n$ and $\Amat(\sigma)=\Amat(\sigma+\Sigma)$ is a $\Sigma$-periodic $n\times n$ matrix \cite{Floquet,Farkas}. The general solution of this system is characterized by the $n\times n$ state transition matrix  $\Xmat(\sigma,s)$ ($s\le\sigma\le s+\Sigma$) since $\zvec_0(\sigma)=\Xmat(\sigma,s)\cvec$ where $\cvec$ is a constant vector. According to Floquet theorem, we have $\Xmat(\sigma,s)=\Mmat(\sigma,s)\exp[\Fmat(\sigma-s)]$ where $\Mmat$ is $n\times n$ and $\Sigma$ periodic with respect to both variables, and $\Fmat$ (of size $n\times n$) is a constant matrix.

Let us now discuss the general linear homogeneous system
\begin{equation}
\totder{\zvec_0}{\sigma}=\Lt\left\{\zvec_0,\sigma\right\}%+\bvec(\sigma)
\label{sifula0}
\end{equation}
where $\Lt$ is a linear operator. There have been some previous attempts to generalize Floquet theory to systems with memory or non-local effects. The most important one is certainly the work of Hale and Verduyn Lunel \cite{Hale1993}. They have derived a generalization of Floquet theory for linear equations of the type \eqref{sifula0} under an important restriction on $\Lt$: that it is $\Sigma$-periodic, i.e., $\Lt\{\zvec,\sigma+\Sigma\}=\Lt\{\zvec,\sigma\}$ for any $\sigma$ and any vector function $\zvec(\sigma)$. However, this set of linear operators does not cover the entire class of systems with memory/non-local effects. Indeed, in many applications the periodicity of the linear operator is guaranteed only if also $\zvec(\sigma)$ is $\Sigma$-periodic. Notice that delay systems, studied in \cite{Simmendinger}, are a particular case of \eqref{sifula0}.

We then demonstrate a theorem that extends Floquet theory to a much larger class of  linear operators $\Lt$. The main assumption is that a positive real $\Sigma$ exists such that $\Lt$ commutes with the translation operator $\T$, defined as the linear operator $\T\left\{ \vvec(\sigma) \right\}=\vvec(\sigma+\Sigma)$. Furthermore, we assume that $\T$ acts on functions $\zvec_0$ defined at least on the interval $[-r+s,s+\Sigma]$ where $s$ is the initial $\sigma$ value for \eqref{sifula0} and $r$ is a real constant.

Denoting as $\mathfrak{L}(\T)$ the set of {\it all} the linear operators $\Lt$ that commute with $\T$, a trivial verification of the definition conditions shows that $\mathfrak{L}(\T)$ is an algebra generated by $\T$  on the topological vector space $\mathcal{H}$, made of functions $\vvec(t):[-r+s,s+\Sigma]\rightarrow\R^n$. Clearly, $\Sigma$-periodic operators are members of this algebra.

%The proof of the theorem relies only on well-known concepts of functional analysis, such as topological vector space and Banach algebra \cite{Rudin2}.

After this brief but necessary introduction we can then formulate the main result of this contribution in the form of the following theorem:
\begin{theorem}[Generalized Floquet theorem] \label{FG}
Let us consider the homogeneous system \eqref{sifula0} with $\Lt\in\mathfrak{L}(\T)$, whose space of solutions is spanned by the state transition matrix $\Xmat(\sigma;s)$ (for $s\le \sigma\le s+\Sigma$). The state transition matrix has then size $n\times p$ with $p\le +\infty$, and can be written as $
%\begin{equation}
\Xmat(\sigma;s)=\Mmat(\sigma;s)\exp[\Fmat(\sigma-s)],
$ %\end{equation}
where $\Mmat(\sigma;s)\in\R^{n\times p}$ is $\Sigma$-periodic in both variables [\idest, $\Mmat(\sigma+\Sigma;s)=\Mmat(\sigma;s)$ and $\Mmat(\sigma;s+\Sigma)=\Mmat(\sigma;s)$] and  $\Fmat\in\R^{p\times p}$ is a constant matrix.
\end{theorem}
\paragraph{Proof.} The first part of the proof is devoted to the definition of the {\it size} of $\Xmat(\sigma;s)$. Let us consider
\begin{equation}
\Lts\left\{\cdot\right\}=\totder{\phantom{\sigma}}{\sigma}\left\{ \cdot\right\}-\Lt\left\{\cdot,\sigma\right\}.
\label{pippo}
\end{equation}
$\Lts$ is a linear operator acting on elements $\zvec_0\in\HT$, the Hilbert space spanned by the functions $\zvec_0:[-r+s,s+\Sigma]\rightarrow\R^n$. Since the functions $\zvec_0$ are defined over the finite interval $[-r+s,s+\Sigma]$, $\HT$ has an infinite countable dimension and according to \cite{Rudin2}, $\ker\left\{\Lts\right\}$ is a countable Hilbert space of dimension $p\le+\infty$. As a consequence $\Xmat(\sigma;s)\in \R^{n\times p}$.

%Because of definition \eqref{pippo} the homogeneous solution $\zvec_0\in\ker\left\{\Lts\right\}$ is a solution of \eqref{sifula}, thus we can build a state transition matrix collecting, as columns of $\Xmat(t;s)$, the elements $\kvec_j(t)$ ($j=1,\dots,p$) of a base of $\ker\{\Lts\}$. {\bf Not necessary: In fact, by definition, the (possibly countably infinite) columns of $\Xmat(t;s)$ span the entire space of solutions of \eqref{sifula}.}

Let us now turn to the {\it form} of $\Xmat(\sigma;s)$. We immediately see that $\Lts\{\Xmat(\sigma;s)\}=\mathbf0$ because it is a state transition matrix for \eqref{sifula0}, and clearly $\T\{\Lts\{\Xmat(\sigma;s)\}\}=\T\{\mathbf0\}=\mathbf0$. Since $\T$ commutes with $\Lt$ and with $\D^\alpha/\D \sigma^\alpha$ for any $\alpha\in\N$, we get
\begin{equation}
\T\{\Lts\{\Xmat(\sigma;s)\}\}=\Lts\{\T\{\Xmat(\sigma;s)\}\}=\mathbf0.
\label{commutation}
\end{equation}

%Turning to the form of $\Xmat(t;s)$, since it is a state transition matrix for \eqref{sifula}
%\begin{equation}
%\totder{\Xmat(t;s)}{t}=\Amat(t)\Xmat(t;s)+\Lt\left\{\Xmat(t;s)\right\}.
%\end{equation}
%The application of the time translation operator, that commutes with $\D/\D t\in\mathfrak{L}(\T)$, yields
%\begin{align}
%\T\left\{ \totder{\Xmat(t;s)}{t}\right\} = \totder{\Xmat(t+T;s)}{t} =  \totder{\T\left\{\Xmat(t;s)\right\}}{t}.
%\label{k1}
%\end{align}

%On the other hand,  a direct verification shows that $\T$ commutes with all the $T$ periodic operators, and by definition of $\Lt\in\mathfrak{L}(\T)$
%\begin{align}
%&\T\left\{ \Amat(t)\Xmat(t;s)+\Lt\left\{\Xmat(t;s)\right\}\right\} \nonumber\\[1ex]
%&\qquad= \Amat(t) \T\left\{\Xmat(t;s)\right\}+\Lt\left\{ \T\left\{\Xmat(t;s)\right\}\right\}.
%\label{k2}
%\end{align}
%Equations \eqref{k1} and \eqref{k2} imply that $\T\left\{\Xmat(t;s)\right\}$ still is a state transition matrix of \eqref{sifula}. Being formed by the basis elements of the kernel of the unbounded operator $\Lts$, the state transition matrix $\T\left\{\Xmat(t;s)\right\}$ can be expressed as a combination of the elements $\kvec_j(t;s)$ of a base of $\ker\{\Lts\}$, \idest\ as a combination of the columns of $\Xmat(t;s)$. The linear operator $\T$ is represented by a matrix, thus such combination is linear. These remarks lead to
Equation \eqref{commutation} implies that $\T\left\{\Xmat(\sigma;s)\right\}$ is a state transition matrix of \eqref{sifula0}. Being formed by the basis elements of the kernel of the linear operator $\Lts$, the state transition matrix $\T\left\{\Xmat(\sigma;s)\right\}$ can be expressed as a combination of the elements $\kvec_j(\sigma)$ of a base of $\ker\{\Lts\}$, \idest\,, as a combination of the columns of $\Xmat(\sigma;s)$. From elementary representation theory, a representation of the linear operator $\T$ is a matrix, therefore such combination is linear. These remarks lead to
\begin{equation}
\T\left\{\Xmat(\sigma;s)\right\}=\Xmat(\sigma+\Sigma;s)=\Xmat(\sigma;s)\Cmat.
\label{k3}
\end{equation}
where the matrix $\Cmat\in\R^{p\times p}$ is constant because $\T$ commutes with $\D^\alpha/\D \sigma^\alpha$.

Since the exponential function is never null, without any loss of generality we can require that the form of $\Xmat(\sigma;s)$ is
$%\begin{equation}
\Xmat(\sigma;s)=\Mmat(\sigma;s)\exp[{\Fmat (\sigma-s)}],
$ %\end{equation}
where $\Fmat\in\R^{p\times p}$ is a constant matrix. Using this in \eqref{k3}
\begin{align}
\T\left\{\Mmat(\sigma;s)\espo{\Fmat (\sigma-s)}\right\}&=\Mmat(\sigma+\Sigma;s)\espo{\Fmat (\sigma+\Sigma-s)}\nonumber\\[1ex]
&=\Mmat(\sigma;s)\espo{\Fmat (\sigma-s)}\Cmat
\label{k4}
\end{align}
that is satisfied only if $\Cmat=\exp({\Fmat \Sigma})$ and $\Mmat(\sigma+\Sigma;s)=\Mmat(\sigma;s)$,
%\begin{equation}
%\Cmat=\espo{\Fmat T}\qquad \up{and} \qquad \Mmat(t+T;s)=\Mmat(t;s)
%\end{equation}
which is the result we set to prove. Repeating the proof with a time translation operator acting on $s$, we finally find that $\Mmat(\sigma;s)$ is periodic also with respect to $s$. This completes the proof of the theorem.
%\end{proof}

Because of the generalized Floquet theorem%~\ref{FG}
, a solution of  \eqref{sifula0} is
\begin{equation}
\zvec_0(\sigma)=\Mmat(\sigma;0)\espo{\Fmat \sigma}\cvec_0'=\Mmat(\sigma;0)\Qmat\espo{\Dmat \sigma}\cvec_0
\label{baraka}
\end{equation}
where $\Dmat\in\R^{p\times p}$ is a diagonal matrix corresponding to the solution of the eigenvalue problem $\Fmat\Qmat=\Qmat\Dmat$, and $\cvec_0'=\Qmat\cvec_0\in\R^m$ is a constant vector. Choosing $\cvec_0=\evec_j$ (the $j$-th element of the canonical unit vector base of $\R^p$), Eq.~\eqref{baraka} yields the general form of the solution
\begin{equation}
\zvec_0(\sigma)=\rvec_j(\sigma)\espo{\lambda_j \sigma}\qquad j=1,\dots,p,
\label{kupenda}
\end{equation}
where $\lambda_j$ is the $j$-th diagonal element of $\Dmat$ (Floquet exponent) and $\rvec_j(\sigma)$ is a $\Sigma$-periodic vector function (Floquet direct eigenvector). Notice that $\Dmat$ is not uniquely defined, since for every $\lambda_{j}$ any complex number $\lambda_{j}+\gei k2\pi/\Sigma$ ($k$ integer) still satisfies \eqref{kupenda} with $\rvec'_j(\sigma)=\rvec_j(\sigma)\exp(-\gei k2\pi\sigma/\Sigma)$. However, the Floquet multipliers $\mu_j=\exp(\lambda_j\Sigma)$ are in number of $p$. Dropping the index $j$, substituting \eqref{kupenda} into \eqref{sifula0} yields
\begin{align}
\totder{\rvec}{\sigma}+\lambda \rvec=\espo{-\lambda \sigma}\Lt\left\{\espo{\lambda \sigma}\rvec,\sigma\right\},
\label{lambda}
\end{align}
where the right hand side is $\Sigma$-periodic because the left hand side is. %Notice that, because of the presence of the term depending on the integro-differential operator $\Lt$, this kind of problems make impossible to use the efficient numerical tehcniques based on multiple shooting for the calculation of the monodromy matrix. In fact, the memory effects due to $\Lt$ make the term dependent on the entire cycle history, \idest\ on the entire set of discretization times in the period.

\paragraph{Dynamical systems with memory.} A first application of the generalized Floquet theorem is to dynamical systems with memory~\cite{Max_review}. We consider nonlinear dynamical systems described by a differential equation of the type
\begin{equation}
\totder{\yvec}{t}=\fvec(\yvec,t)+\Ltp\left\{\gvec(\yvec,t),t\right\}
\label{nonlinsys}
\end{equation}
where $\yvec: \R\rightarrow \R^n$, $\fvec,\gvec: \R^{n+1}\rightarrow \R^n$ are nonlinear vector fields with $\fvec$ $T$-periodic with respect to $t$, and $\Ltp$ is a generic, linear integro-differential operator such that $\Ltp\left\{\gvec(\yvec,t),t\right\}$ is $T$-periodic if $\yvec(t)$ is $T$-periodic. We assume that \eqref{nonlinsys} admits a $T$-periodic solution $\yvec_\up{S}$ (limit cycle).

The study of the perturbation of the limit cycle $\yvec_\up{S}$ is important both for assessing the stability of the solution and for analyzing the effect of small-amplitude forcing terms such as the Langevin sources used to represent fluctuation effects. Perturbations are studied by linearizing \eqref{nonlinsys} around $\yvec_\up{S}$ yielding a linear periodic time-varying system
\begin{equation}
\totder{\zvec}{t}=\Amat(t)\zvec+\Ltp\left\{\Bmat(t)\zvec,t\right\}+\bvec(t)
\label{sifula}
\end{equation}
where $\zvec=\yvec-\yvec_\up{S}: \R\rightarrow\R^n$, $\Amat(t)$ is the $T$-periodic Jacobian matrix of $\fvec(\cdot,t)$ and $\Bmat(t)$ is the Jacobian matrix of $\gvec(\cdot,t)$ both evaluated in the limit cycle, and $\bvec(t)$ is the perturbation forcing term. Finally, $\Lt$ is the {linear operator} defined as $
%\begin{equation}
\Lt\left\{\cdot,t\right\}=\Amat(t)\cdot + \Ltp\left\{\Bmat(t)\cdot\right\}
%\end{equation}
$. The solution of the homogeneous equation associated to \eqref{sifula} (i.e., the case $\bvec(t)=\zerovec$) is denoted as $\zvec_0(t)$.

Dynamical systems with memory are included in \eqref{sifula} as convolution operators:
\begin{equation}
\Ltp\left\{\Bmat\zvec,t\right\}=\int_{t-r}^t \Kmat(t;\tau)\zvec(\tau)~\D\tau
\label{memory}
\end{equation}
where $r>0$ and $\Kmat$ is a matrix such that the integral of $\Kmat(t;\tau)$ with respect to $\tau$ is bounded for every $t$, and $\Kmat(t+T;\tau+T)=\Kmat(t;\tau)$. The latter condition (sometimes called bi-periodicity) is not enough to guarantee the periodicity of $\Lt$ for any $\zvec$ (see above), but it is just true when $\zvec$ is $T$-periodic. Furthermore, it is simple to prove that $\Lt$ commutes with $\T$, thus implying that the generalized Floquet theorem applies for the solution of the homogeneous part of \eqref{sifula}.

The asymptotic stability analysis of $\yvec_\up{S}(t)$, therefore, clearly depends on the computation of all the eigenvalues of $\Fmat$: this procedure depends on the explicit form of $\Lt$, and in general leads to $p$ classes of Floquet exponents $\lambda$~\footnote{Notice that according to the generalized Floquet theorem $p$ may also be infinite.}. According to the previous discussion, each class is characterized by an infinite set of complex numbers having the same real part, and imaginary parts differing by integer multiples of $2\pi/T$.

Another important remark is that direct substitution proves that the general solution of \eqref{sifula} can be written as
\begin{equation}
\zvec(\sigma)=\Xmat(\sigma;0)\cvec_0+\int_0^\sigma \Xmat(\sigma;\eta)\bvec'(\eta)~\D\eta
\end{equation}
where $\cvec_0$ is a constant vector of size $p$, and $\bvec'(\sigma)$ is a solution of
\begin{align}
\Xmat(0;0)\bvec'(\sigma)&=\bvec(\sigma)+\Lt\left\{\int_0^\sigma \Xmat(\sigma;\eta)\bvec'(\eta)~\D\eta\right\} \nonumber\\
&\qquad- \int_0^\sigma \Lt\left\{\Xmat(\sigma;\eta)\right\}\bvec'(\eta)~\D\eta.
\end{align}
This result can be used, for example, to study phase and amplitude noise in autonomous (self-oscillating) systems, extending \cite{noiFN,noi}, e.g., to the case of electronic oscillators including distributed circuit elements (namely, transmission lines).

Finally, the results of the theorem can be extended to operators of type \eqref{memory} for  $r\rightarrow+\infty$ under the following conditions: when $r\rightarrow+\infty$, because of the properties of $\Kmat$, $\Kmat(\sigma,\eta)\rightarrow 0$ for $\eta\rightarrow -\infty$ and any $\sigma$. For any $\epsilon>0$ there is an $r=r(\epsilon)$ such that $\abs{\int_{-\infty}^{\sigma-r} \Kmat(\sigma;\eta) \yvec_\up{S}~\D\eta}<\epsilon$ being $\yvec_\up{S}$ $\Sigma$-periodic and integrable. Thus, for any $\epsilon$ we can approximate an infinite memory system with \eqref{memory}. In this sense our theorem guarantees that the solution of the homogeneous part of system \eqref{sifula} is given by \eqref{kupenda} for $\sigma\in[-r+s,s+\Sigma]$ and for any $r<+\infty$.

\paragraph{Quantum systems with non-local potentials.} The second application we consider relates to single-particle quantum systems (with mass $m$) characterized by the 3D time-independent Schrodinger equation where the potential energy linear operator $\V\{\psi,\xvec\}$  commutes with a properly defined 3D extension of the translation operator $\T$ ($\psi(\xvec)$ is the particle wave-function).

We assume that three linearly independent spatial vectors $\nvec_\alpha$ ($\alpha=1,2,3$) exist such that $\avec=k_1\nvec_1+k_2\nvec_2+ k_3\nvec_3$ ($k_\alpha\in\mathbb{Z}$) defines the direct space lattice in which the particle is embedded. The generalized 3D translation operator is defined as $\Tb_{\nvec_\alpha}=\T_{n_{1,\alpha}}\T_{n_{2,\alpha}}\T_{n_{3,\alpha}}$, where $\T_{n_{i,\alpha}}$ is the translation operator acting only on the $i$-th ($i=1,2,3$) component of $\nvec_\alpha$~\footnote{Clearly the three factors of the composition product commute with each other.}.

In order to apply the generalized Floquet theorem, we need to reduce the problem from 3D to 1D. The first step is to define a transformation matrix $\Rmat$ that transforms the three basis vectors $\nvec_\alpha$ into the canonical basis of $\mathbb{R}^3$: $\Rmat\nvec_\alpha=\evec_\alpha$. Correspondingly, we have the transformed spatial variable $\tilde{\xvec}=\Rmat\xvec$ and the transformed second order operator defining the differential part of Schroedinger equation $\nabla^2=(\Rmat\tilde\nabla)^2$, where the power should be intended as a scalar product. Also the shift operator is modified by the variable transformation: in the canonical base the shift operator becomes  $\Tb_{\evec_\alpha}=\T_{\tilde{x}_\alpha}$ ($\alpha=1,2,3$) where the period $\Sigma$ is equal to 1 in all the three spatial directions. Of course $\Vtilde\{\tilde\psi,\tilde{\xvec}\}=\V\{\psi,\xvec\}$ commutes with $\Tb_{\evec_\alpha}$ ($\alpha=1,2,3$), and $\tilde\psi(\tilde\xvec)=\psi(\Rmat^{-1}\tilde\xvec)=\psi(\xvec)$.

The orthogonalization of the direct lattice made possible by $\Rmat$ allows us to transform the 3D equation into three  different couples of 1D first order differential equations. Each of the couples is expressed as differential with respect to one transformed spatial variable at a time, but of course the solution must be the same for all the possible choices. For instance, with respect to $\tilde{x}_1$ we have ($\Rvec_i$ is the $i$-th column of $\Rmat$, $\tilde\partial_j$ the partial derivative with respect to $\tilde x_j$)
\begin{subequations}
\label{puffo}
\begin{align}
&\tilde\partial_1\tilde\psi=\tilde\chi_1
\label{puffo1}\\
&\tilde\partial_1\tilde\chi_1=\fracd{1}{\Rvec_1^2}\fracd{2m}{\hbar^2}\left\{ \Vtilde\{\tilde{\psi}\}- \left[E+\fracd{\hbar^2}{2m}\left( \Rvec_2\tilde\partial_2+ \Rvec_3\tilde\partial_3 \right)^2\right]\tilde\psi\right\} \nonumber\\
&\quad
+\fracd{2}{\Rvec_1^2} \Rvec_1\cdot\left( \Rvec_2\tilde\partial_2
+ \Rvec_3\tilde\partial_3 \right)\tilde\chi_1
\label{puffo2}
\end{align}
\end{subequations}
where $E$ is the particle energy. Similar systems can be derived with reference to spatial derivatives with respect to $\tilde{x}_2$ and $\tilde{x}_3$, and introducing $\tilde\chi_\beta=\tilde\partial_\beta\tilde\psi$ ($\beta=2,3$).

Since the right hand side of \eqref{puffo} commutes with $\T_{\tilde{x}_1}$ for $\Sigma=1$, the generalized Floquet theorem can be applied to obtain the general solution
\begin{subequations}
\label{p}
\begin{equation}
\tilde\psi(\tilde\xvec)=c_1(\tilde x_2,\tilde x_3) \espo{\lambda_1(E)\tilde x_1} \tilde u(\tilde\xvec;E)
\label{p1}
\end{equation}
where $\tilde u(\tilde \xvec;E)$ is periodic with period equal to 1 as a function of $\tilde x_1$. Notice that $\lambda_1$ depends on $E$ only (and not on $x_2$ and $x_3$) because $\exp(\lambda_1)$ pertains to the spectrum of $\T_{\tilde{x}_1}$.  Similarly, by solving the other two equations
\begin{align}
\tilde\psi(\tilde\xvec)&=c_2(\tilde x_1,\tilde x_3) \espo{\lambda_2(E)\tilde x_2} \tilde u(\tilde\xvec;E) \label{p2} \\
\tilde\psi(\tilde\xvec)&=c_3(\tilde x_1,\tilde x_2) \espo{\lambda_3(E)\tilde x_3} \tilde u(\tilde\xvec;E).
\label{p3}
\end{align}
\end{subequations}
Imposing that \eqref{p} defines in three different ways the same wavefunction, it can be easily shown, expressing the $c$ coefficients in exponential form and equating \eqref{p1}--\eqref{p3}, that the three coefficients $c_i$ are equal to $\exp(\lambda_j\tilde x_j)\exp(\lambda_k\tilde x_k)$ (with $j\neq k\neq i$). Therefore, collecting the $\lambda_j(E)$ in a vector $\lambdavec$ of size 3, we have the general solution of the 3D Schroedinger equation in real space as %($u(\xvec;E)=\tilde{u}(\Rmat\xvec;E)$)
\begin{equation}
\psi(\xvec)=\espo{\lambdavec(E)\cdot\Rmat\xvec}u(\xvec;E)=\espo{\Rmat^\up{T}\lambdavec(E)\cdot\xvec}u(\xvec;E)
\end{equation}

Since \eqref{puffo1} and \eqref{puffo2} are of size 2, there are at least two vector Floquet exponents for each $E$ value. However, for a non-local potential energy $\V$, according to the generalized Floquet theorem, for each value of $E$ more values (possibly infinite) than the pure local potential case of $\lambdavec(E)$ are possible depending on the specific shape of $\V$.

Imposing that the normalization constant over a unit cell of the direct lattice is independent of the chosen cell (invariance with respect to a combination with integer coefficients of the $\nvec_j$ unit vectors), one finds $\Rmat^\up{T}\lambdavec(E)=\pm\gei\kvec(E)$ that reduces to progressive and regressive plane waves for a local potential energy ($p=2$), $\kvec$ being the Bloch function wavevector (and, thus, $\kvec(E)$ the corresponding band diagram). In other words, the dispersion relation $E(\kvec)$ is an even function. Notice that, e.g., in the presence of a magnetic field, time-reversal invariance is broken thus lifting Kramers degeneracy. A magnetic field implies the appearance of an imaginary term in the Hamiltonian that, in turn, no longer yields conjugate imaginary $\lambdavec(E)$ values (and the following parity of the dispersion relation as a function of the Bloch wavevector).

Thus, we can summarize the two main results we find from including the nonlocal potential into the Schroedinger equation as follows: i) we still have Bloch waves as solutions of the time-independent Schroedinger equation, ii) the dispersion relation $E(\kvec)$ can qualitatively change according to the number of $\kvec$ vectors associated to the energy.
The latter result is too complex to be fully described here, so to provide physical insight we limit the description to the 1D case, although the same conclusions can be extended to the 3D case. For a local potential in the absence of magnetic field, since $E(k)$ is an even function of $k$, and for each $k$ value only two allowed energy states are available (for $\pm k$) with the same energy, the band diagram can exhibit an extremum only at the center of the first Brillouin zone (FBZ) or at its borders: the band diagram is monotonic as a function of ${k}$ in each half of the FBZ.

On the other hand, for non local $\V$ we may have $p>2$ -- although always an even number, because the normalization condition involves the presence of imaginary conjugate $\lambda$ values -- thus leading to the possible occurrence of non-monotonic dispersion relations $E(k)$ in each half of the FBZ (i.e., the energy bands may show extrema also within the half FBZ, while this is forbidden for local potentials). Furthermore, by tuning the strength/shape of the non local part of the potential energy, energy band crossings in the first Brilluoin zone may be observed, leading to
possible quantum phase transitions.

%Combining the three cases, and taking into account the physical conditions on the normalization of the eigenfunctions, the resulting general solution of the Schroedinger equation in the original coordinates is
%\begin{equation}
%\psi(\xvec)=\espo{\pm\gei\kvec(E)\cdot\xvec}\uvec(\xvec,E)
%\end{equation}
%where $\uvec$ has the same periodicity of $\Vtilde(\xvec)$ and $\kvec(E)$ is a set of real 3D vectors whose size depends on the specific shape of $\Vtilde(\rvec)$: for a nonlocal potential, two vectors are found with opposite sign. However, for nonlocal potentials such number may become larger than 2.

\paragraph{Conclusions.} In conclusion, we have stated and proved a theorem that generalizes the standard Floquet theory
to a wide class of systems, that includes those whose state depends on the full past dynamics, and those whose  potential energy profile, albeit spatially periodic, is non local. Such a theorem greatly expands Floquet analysis to systems that
were previously inaccessible by such a theory, making it a tool with even greater potential in several areas of science and technology.

For instance, the extension to dynamical systems with memory characterized by a time periodic limit cycle opens the way to a consistent stability theory for such systems. Examples of applications are related to, e.g., all those electronic circuits containing distributed elements, such as transmission lines. In this area, Floquet theory is also exploited to assess the properties of autonomous systems including the effect of stochastic fluctuations, leading to the phase and amplitude noise analysis of a wide range of oscillators, such as electronic periodic signal generators \cite{noi} or optical sources such as lasers \cite{lasernoise,Demir1}.
Furthermore, this work opens the way to the analysis of Bloch states in condensed matter systems characterized by non-local potentials. Interesting results are already available even at this early stage of application, such as the possible occurrence of multi-valued dispersion relations and the possibility to tune the non-locality to
 induce quantum phase transitions.

%\clearpage

{\it Acknowledgments --} This work has been partially supported by the Spanish government through project TEC2011-14253-E, by the NSF grant No. DMR-0802830, and the Center for Magnetic Recording Research at UCSD.
%\section{Example}

%\appendix
%\section{}

%\begin{acknowledgments}
%\end{acknowledgments}

%\clearpage

\bibliography{IEEEabrv,APStraversa1}

\end{document}